\documentclass[a4paper,12pt]{article}
\usepackage[centertags]{amsmath}
\usepackage[french,english]{babel}
\usepackage{amsfonts}
\usepackage{amsmath,amssymb,graphics,amscd}
\usepackage{amsthm}
\usepackage{color}
\definecolor{mycolor}{rgb}{1, 0, 0}
\newcounter{jolist}

\newcommand{\semi}{{\,\rule[.1pt]{.4pt}{5.3pt}\hskip-1.9pt\times}}

\newcommand{\al}{\alpha}

\newcommand{\eps}{\epsilon}
\newcommand{\om}{\omega}

\newcommand{\lam}{{\lambda}}
\newcommand{\Lam}{{\Lambda}}

\newcommand{\Lhg}{\widehat{\mathfrak{g}}}
\newcommand{\Lhh}{\widehat{\mathfrak{h}}}
\newcommand{\Lhb}{\widehat{\mathfrak{b}}}
\newcommand{\wDelta}{\widehat{{\Delta}}}
\newcommand{\wPhi}{\widehat{{\Phi}}}
\newcommand{\ph}{\widehat{P}}
\newcommand{\charc}{\hbox{\rm Char}\,}
\newcommand{\waff}{{W^{\rm aff}}}
\newcommand{\wiff}{{\widetilde{W}^{\rm aff}}}

\newcommand{\semc}{{\rm sc}}

\newcommand{\bn}{{\mathbb N}}
\newcommand{\br}{{\mathbb R}}
\newcommand{\bz}{{\mathbb Z}}

\newcommand{\bc}{{\mathbb C}}

\newcommand{\La}{{\mathfrak{a}}}
\newcommand{\Lg}{{\mathfrak{g}}}
\newcommand{\Lb}{{\mathfrak{b}}}
\newcommand{\Lh}{{\mathfrak{h}}}

\newcommand{\Ld}{{\mathfrak{d}}}

\newtheorem{thm}{Theorem}
\newtheorem{dfn}{Definition}

\newtheorem{conj}{Conjecture}
\newtheorem{rem}{Remark}

\newtheorem{lem}{Lemma}
\newtheorem{coro}{Corollary}

\def\proof{\noindent{\it Proof. }}
\def\endpf{\hfill$\bullet$\vskip 5pt}

\begin{document}
\title{Tensor product structure of affine Demazure modules and limit constructions}
\author{G. Fourier$^*$ and P. Littelmann\footnote{This research has been partially
supported by the EC TMR network "LieGrits", contract
MRTN-CT 2003-505078. {\sl 2000 Mathematics Subject Classification. 22E46, 14M15.}}}
\maketitle
\bigskip
\begin{abstract}
Let $\Lg$ be a simple complex Lie algebra, we denote by $\Lhg$ the
affine Kac--Moody algebra associated to the extended Dynkin
diagram of $\Lg$. Let $\Lambda_0$ be the fundamental weight of
$\Lhg$ corresponding to the additional node of the extended Dynkin
diagram. For a dominant integral $\Lg$--coweight $\lam^\vee$, the Demazure submodule
$V_{-\lam^\vee}(m\Lam_0)$ is a $\Lg$--module. We provide a
description of the $\Lg$--module structure
as a tensor product of ``smaller'' Demazure modules. More precisely,
for any partition of $\lam^\vee=\sum_j \lam_j^\vee$ as a sum of dominant
integral $\Lg$--coweights,  the Demazure module is (as $\Lg$--module)
isomorphic to
$\bigotimes_j V_{-\lam^\vee_j}(m\Lam_0)$. For the ``smallest'' case,
$\lam^\vee=\om^\vee$ a fundamental coweight, we provide for $\Lg$ of
classical type a decomposition of $V_{-\om^\vee}(m\Lam_0)$ into
irreducible $\Lg$--modules, so this can be viewed as a natural generalization
of the decomposition formulas in \cite{KMOTU} and \cite{Magyar}.
A comparison with the $U_q(\Lg)$--characters of certain finite
dimensional $U_q'(\Lhg)$--modules (Kirillov--Reshetikhin--modules)
suggests  furthermore that all quantized Demazure modules
$V_{-\lam^\vee,q}(m\Lam_0)$ can be naturally endowed with the
structure of a $U_q'(\Lhg)$--module. We prove, in the classical case (and for a lot of non-classical cases),
a conjecture by Kashiwara \cite{Kashiwara2}, that the ``smallest'' Demazure modules are,
when viewed as $\Lg$-modules, isomorphic to some KR-modules.
For an integral dominant $\Lhg$--weight $\Lambda$ let
$V(\Lam)$ be the corresponding irreducible $\Lhg$--representation.
Using the tensor product decomposition for Demazure modules,
we give a description of the $\Lg$--module structure of $V(\Lam)$
as a semi-infinite tensor product of finite dimensional $\Lg$--modules.
The case of twisted affine Kac-Moody algebras can be treated in the same way, some
details are worked out in the last section.
\end{abstract}
\section*{Introduction}
Let $\Lg$ be a simple complex Lie algebra, we denote by $\Lhg$ the
affine Kac--Moody algebra associated to the extended Dynkin
diagram of $\Lg$. (The twisted case is considered separately in the last section).
Let $\Lambda_0$ be the fundamental weight of
$\Lhg$ corresponding to the additional node of the extended Dynkin
diagram. The basic representation $V(\Lam_0)$ is one of the most
important representations of $\Lhg$ because its structure
determines strongly the structure of all other highest weight
representations $V(\Lam)$, $\Lam$ an arbitrary dominant integral
weight for $\Lhg$.

Let $P^\vee$ be the coweight lattice of $\Lg$. An element
$\lam^{\vee}$ in the coroot lattice, can be viewed as an element
of the affine Weyl group $\waff$ (see section 1), and one can associate to
$\lam^\vee$ the Demazure submodule $V_{\lam^\vee}(\Lam)$ of
$V(\Lam)$ (see section 2).

Actually, this construction generalizes to arbitrary $\lam^\vee\in
P^\vee$ in the following way: one can write $\lam^\vee$ as
$w\sigma\in \wiff$ in the extended affine Weyl group, where $w\in
\waff$ and $\sigma$ corresponds to an automorphism of the Dynkin
diagram of $\Lhg$. Denote  by $V_{\lam^\vee}(\Lam)$ the Demazure
submodule $V_w(\Lam')$ of the highest weight module $V(\Lam')$,
where $\Lam'= \sigma(\Lam)$.

If $\lam^\vee$ is a dominant coweight, then the Demazure module
$V_{-\lam^\vee}(m\Lam_0)$ is in fact a $\Lg$--module, and it is
interesting to study its structure as $\Lg$--module. So one would
like to get a restriction formula expressing $V_{-\lam^\vee}(m\Lam_0)$
as a direct sum of simple $\Lg$--representations. We write
$\overline{V}_{-\lam^\vee}(m\Lam_0)$ for the Demazure module viewed
as a $\Lg$--module.

A first reduction step is the following theorem describing
the Demazure module as a tensor product.
Such a decomposition formula for Demazure modules was first
observed by Sanderson~\cite{Sanderson} in the affine rank two case,
and was later studied in the case of classical groups in the framework
of perfect crystals for example in \cite{KMOTU2},  see \cite{HongKang}
for a more complete account.
We provide in this article a description of the Demazure module
as a tensor product of modules of the same type,
but for ``smaller coweights''. More precisely, let $\lam^\vee$ be a
dominant coweight and suppose we are given a decomposition
$$
\lam^\vee=\lam_1^\vee+\lam_2^\vee+\ldots+\lam_r^\vee
$$
of $\lam^\vee$ as a sum of dominant coweights. The following
theorem is a generalization of a result in \cite{Magyar}, where
the statement has been proved in the case $m=1$ and under the additional
assumption that all the $\lam_i^\vee$ are minuscule fundamental
weights, and the decomposition formulas in
\cite{HKOTT}, \cite{HKKOTY}, \cite{KMOTU} and \cite{KMOTU2},
where in the framework of perfect crystals for classical groups
many cases have been discussed. (The corresponding version
for a twisted Kac-Moody algebra can be found in section 4.)
\vskip 8pt\noindent
{\bf Theorem 1.} {\it
For all $m\ge 1$, we have an isomorphism of $\Lg$--representations
between the Demazure module $\overline{V}_{{-\lam^\vee}}(m\Lam_0)$
and the tensor product of Demazure modules:}
$$
\overline{V}_{{-\lam^\vee}}(m\Lam_0)\simeq
\overline{V}_{{-\lam_1^\vee}}(m\Lam_0)\otimes
\overline{V}_{{-\lam_2^\vee}}(m\Lam_0)
\otimes\cdots\otimes
\overline{V}_{{-\lam_r^\vee}}(m\Lam_0).
$$
\vskip 8pt
Of course, to analyse the structure of
$\overline{V}_{{-\lam^\vee}}(m\Lam_0)$ as a $\Lg$--module, the
simplest way is to take a decomposition of $\lam^\vee$ as a sum of
fundamental coweights $\lam^\vee=\sum a_i\om_i^\vee$. So by
Theorem 1, it remains to describe the structure of the
$\overline{V}_{{-\om_i^\vee}}(m\Lam_0)$ as a $\Lg$--module. We
give such a description below for all fundamental coweights for the
classical groups. For the exceptional groups we give the
decomposition in the cases interesting for the limit constructions
considered later. The enumeration of the fundamental weights is as
in \cite{Bourbaki}, we write $\om_0$ for the trivial weight. For more
details on the notation see section 2, we only recall here
that we use the abbreviations $V_{-\om^\vee}(m\Lam_0)$ and
$V_{-\nu(\om^\vee)}(m\Lam_0)$ for the Demazure submodule
associated to the translation $t_{-\nu(\om^\vee)}$,
viewed as an element in the extended affine Weyl group.
We write $V^*$ for the contragradient dual of a representation $V$.
Many special cases of the list below have been calculated before,
for example by the Kyoto school (\cite{Kashiwara2}, \cite{KMOTU}, \cite{KMOTU2}, \cite{Yamane}).
(The corresponding version
for a twisted Kac-Moody algebra can be found in section 4.)
\vskip 5pt\noindent
{\bf Theorem 2.} {\it Let $\omega^\vee$ be a fundamental coweight and
let  ${V}_{-\omega^\vee}(m\Lam_0)$ be the associated Demazure
module. Viewed as a $\Lg$--module, $\overline{V}_{-\omega^\vee}(m\Lam_0)$
decomposes into the direct sum of irreducible $\Lg$--modules as follows:
\begin{itemize}
\item Type ${\tt A}_n$: $\overline{V}_{{-\om_i^\vee}}(m\Lam_0)
=\overline{V}_{{-\om_i}}(m\Lam_0)
\simeq V(m\omega_i)^*$ as ${\mathfrak{sl}}_n$--module for all
$i=1,\ldots,n$.
\item Type ${\tt B}_n$:  Set $\theta=0$ for $i$ even and $\theta=1$ for $i$ odd, then
we have for $1\le i<n$:
$$
\overline{V}_{-\om_i^\vee}(m\Lam_0)=
\overline{V}_{-\om_i}(m\Lam_0) \simeq
\bigoplus_{a_{i}+a_{i-2}+\ldots + a_{\theta}= m}
V(a_i\omega_i+a_{i-2}\om_{i-2}+\ldots+a_{\theta}\omega_{\theta})
$$
and for $i=n$:
$$
\overline{V}_{-\om_n^\vee}(m\Lam_0)=
\overline{V}_{-2\om_n}(m\Lam_0) \simeq
\bigoplus_{a_{n}+a_{n-2}+\ldots + a_{\theta}= m}
V(2a_n\omega_n+a_{n-2}\om_{n-2}+\ldots+a_{\theta}\omega_{\theta})
$$
as ${\mathfrak{so}}_{2n+1}$--module.
\item Type ${\tt C}_n$: for $j<n$ we have
$$
\overline{V}_{{-\om_j^\vee}}(m\Lam_0)=\overline{V}_{{-2\om_j }}(m\Lam_0)
\simeq\bigoplus_{a_1+\ldots+a_j\le m}{V}(2a_1\om_1+\ldots +2a_j\om_j)
$$
and for $j=n$:
$\overline{V}_{{-\om_n^\vee}}(m\Lam_0)=\overline{V}_{{-\om_n}}(m\Lam_0)
\simeq V(m\omega_n)$
as ${\mathfrak{sp}}_n$--module.
\item Type ${\tt D}_n$:
Set $\theta=0$ for $i$ even and $\theta=1$ for $i$ odd, then
we have for $2\le i\le n-2$:
$$
\overline{V}_{-\om_{i}^\vee}(m\Lam_0)=
\overline{V}_{-\om_{i}}(m\Lam_0)\simeq \bigoplus_{a_i+a_{i-2}+\ldots+a_{\theta}=m}
V(a_i\om_i+a_{i-2}\om_{i-2}+\ldots+a_{\theta}\om_{\theta})
$$
and for $i=1,n-1,n$:
$
\overline{V}_{{-\om_i^\vee}}(m\Lam_0) =\overline{V}_{{-\om_i}}(m\Lam_0)
\simeq V(m\omega_i)^*
$
as ${\mathfrak{so}}_{2n}$--module.
\item Type ${\tt E}_6$:
\begin{itemize}
\item[] $\overline{V}_{{-\om_i^\vee}}(m\Lam_0)=\overline{V}_{{-\om_i}}(m\Lam_0)
\simeq V(m\omega_i)^*$ for $i=1,6$%
\item[] $\overline{V}_{{-\om_2^\vee}}(m\Lam_0)=\overline{V}_{{-\om_2}}(m\Lam_0)
\simeq \bigoplus\limits_{r = 0}^{m} V (r \om_2)$ %
\end{itemize}
as ${\tt E}_6$--module.
\item Type ${\tt E}_7$:
\begin{itemize}
\item[] $\overline{V}_{{-\om_7^\vee}}(m\Lam_0)=\overline{V}_{{-\om_7}}(m\Lam_0)
\simeq V(m\omega_7)$%
\item[] $\overline{V}_{{-\om_1^\vee}}(m\Lam_0)=\overline{V}_{{-\om_1}}(m\Lam_0)
\simeq \bigoplus\limits_{r = 0}^{m} V (r \om_1)$ %
\end{itemize}
as ${\tt E}_7$--module.
\item Type ${\tt E}_8$: $\overline{V}_{{-\om_8^\vee}}(m\Lam_0)=\overline{V}_{{-\om_8}}(m\Lam_0)
\simeq \bigoplus\limits_{r = 0}^{m} V (r \om_8)$ %
as ${\tt E}_8$--module.
\item Type ${\tt F}_4$:
$$
\overline{V}_{{-\om_1^\vee}}(m\Lam_0)=
\overline{V}_{{-\om_1}}(m\Lam_0)
\simeq \bigoplus\limits_{r = 0}^{m} V (r \om_1)
$$
and %
$$
\overline{V}_{{-\om_4^\vee}}(m\Lam_0)=
\overline{V}_{{-2\om_4}}(m\Lam_0)
\simeq \bigoplus\limits_{r +s \le m} V (r \om_1+s2\om_4)
$$
as ${\tt F}_4$--module.
\item Type ${\tt G}_2$:
$\overline{V}_{{-\om_2^\vee}}(m\Lam_0)=\overline{V}_{{-\om_2}}(m\Lam_0)
\simeq \bigoplus\limits_{r = 0}^{m} V (r \om_2)$ %
as ${\tt G}_2$--module.
\end{itemize}
}
There is a very interesting conjectural connection
with certain $U'_q (\Lhg)$--modules. Here $U'_q (\Lhg)$ denotes the quantized
affine algebra without derivation. \\
Let $KR(m\om_i)$ be the Kirillov--Reshetikhin--module for a multiple of a fundamental weight of $\Lg$, for the precise definition see \cite{HKOTT}, it is irreducible as $U'_q (\Lhg)$ and the highest weight, when viewed as a $U_q(\Lg)$-module, is $m \om_i$. In \cite{Kashiwara} Kashiwara introduced the notion of a good $U'_q (\Lhg)$--module,
 which, roughly speaking, is an irreducible finite
dimensional $U'_q (\Lhg)$--module
with a crystal basis and a global basis, and he proved that the tensor product
of good modules is a good module. It is conjectured that the KR-modules are good. For all fundamental $\Lg$--weights $\om_i$
Kashiwara constructed this irreducible finite-dimensional integrable $U'_q (\Lhg)$--module
$KR(\om_i)$ and showed that it is good and even more that the crystal
is isomorphic to a certain generalized Demazure crystal as a $\Lg$--crystal.

Let $c_{k}^{\vee} = \frac{a_k}{a_{k}^{\vee}}$
(for the definition of the $a_k$ see section 1.1) and $l \in \mathbb{N}$.
Let $KR(l c_k^{\vee} \om_k)$ be the Kirillov--Reshetikhin--module
for $U_{q}'(\Lhg)$ associated to the weight
$l c_{k}^{\vee} \om_{k}$. It is more generally conjectured that
the $KR(l c_k^{\vee} \om_k)$ the crystal is isomorphic to the crystal of a Demazure module,
after omitting the $0$--arrows in both crystals \cite{Kashiwara2}.
Chari and the Kyoto school have calculated for
classical Lie-algebras and some fundamental weights for non-classical
Lie-algebras the decomposition of the Kirillov--Reshetikhin module
$KR(lc_{k}^{\vee}\om_k))$ into irreducible $U_q(\Lg)$--modules \cite{Chari}. By comparing
the $U_q(\Lg)$--structure of the Kirillov--Reshetikhin module
$KR(lc_{k}^{\vee}\om_k))$ with the list in Theorem~2 we conclude:
\begin{coro}
In all cases stated in Theorem~2,
the Demazure module $(\overline{V}_{{-\om^\vee},q}(l\Lam_0))$
and the Kirillov--Reshetikhin module $KR( l c_{k}^{\vee} {\om_k}^*)$
are, as $U_q(\Lg)$--modules, isomorphic.
\end{coro}
In particular, if $KR( l c_{k}^{\vee} {\om_k}^*)$ has a crystal basis,
then the crystal is isomomorphic to the crystal of
$(V_{-\om_k^{\vee}}( l \Lam_0))$ after
omitting the arrows with label zero.
By using the $U_q(\Lg)$--module isomorphism, we see that
(in the cases above) the quantized Demazure modules $V_{-\om_k^{\vee},q}( l \Lam_0)$
can be equipped with the structure of an irreducible $U'_q (\Lhg)$--module. In fact, using the
Theorem~1, we see that for classical groups all
quantized Demazure modules $V_{-\lam^{\vee},q}( l \Lam_0)$, $\lam^\vee$ a dominant coweight,
can be equipped with the structure of an $U'_q (\Lhg)$--module. Of course, in the
exceptional case the same argument shows that when $\lam^\vee$ can be written
as linear combinations of the fundamental weights listed in Theorem~2, then
again $V_{-\lam^{\vee},q}( l \Lam_0)$ can be equipped with the structure of an
$U'_q (\Lhg)$--module. This leads to the following:
\begin{conj}\rm
Let $\Lg$ be a semisimple Lie algebra, let $U_q (\Lhg)$
be the associated untwisted quantum affine algebra and let
$U'_q(\Lhg)$ be its subalgebra without derivation. For all
dominant coweights $\lam^\vee$ and for all $l > 0$, the Demazure module
$V_{-\lam^{\vee},q}( l \Lam_0)$ can be endowed with the structure
of a $U'_q(\Lhg)$--module admitting a crystal basis. Its crystal graph
is isomorphic to the crystal of the Demazure module, after omitting the
arrows labelled with zero.
\end{conj}

The tensor decomposition structure in Theorem~1
holds in the following more general situation. Let $\Lam_i$,
$1\le i\le n$, be a fundamental weight of $\Lhg$ such that
the corresponding coweight $\om_i^\vee$ is minuscule.
Let $\lam^\vee$ be a dominant coweight and suppose we
are given a decomposition
$$
\lam^\vee=\om^\vee_i+\lam_2^\vee+\ldots+\lam_r^\vee
$$
of $\lam^\vee$ as a sum of dominant coweights and denote $\om_i^*$
the highest weight of the irreducible $\Lg$--module $V(\om_i)^*$.
\vskip 8pt\noindent
{\bf Theorem 1 A.} {\it For all $m\ge ~0$ and
$s\ge 1$, we have an isomorphism of $\Lg$--representations between
the Demazure module $\overline{V}_{{-\lam^\vee}}(m\Lam_0+s\Lam_i)$
and the tensor product of Demazure modules:}
$$
\overline{V}_{{-\lam^\vee}}(m\Lam_0+s\Lam_i)\simeq
{V}(s\om_i^*)\otimes
\overline{V}_{{-\lam_2^\vee}}((m+s)\Lam_0)
\otimes\cdots\otimes
\overline{V}_{{-\lam_r^\vee}}((m+s)\Lam_0).
$$
\vskip 8pt\noindent
Let $\Lam$ be an arbitrary dominant integral weight for $\Lhg$.
The $\Lhg$--module $V(\Lambda)$ is the direct limit of the
Demazure-modules $V_{- N \lambda^{\vee}}(\Lambda)$ for some
dominant, integral, nonzero coweight of $\Lg$. We give a construction of
the $\Lg$--module $\overline{V}(\Lambda)$ as a direct limit of
tensor products of Demazure modules. This has been done before
in the case of classical Lie-algebras for $\Lam=r \Lam_0$ (and
corresponding weights obtained by automorphisms as in the statement
of Theorem 2) by Kang, Kashiwara, Kuniba, Misra et al. \cite{HongKang}, \cite{KMOTU} via
the theory of perfect crystals. In addition they have also
considered some special weights in the case of non-classical groups.
For ${\tt G}_2$, such a construction has been given by Yamane
\cite{Yamane}. For the Lie algebras of type ${\tt E}_6$ and ${\tt
E}_7$ a construction (only for the case $\Lam=\Lam_0$) was given
by Peter Magyar \cite{Magyar} using the path model.

We provide in this article such a direct limit construction for arbitrary
simple Lie algebras $\Lg$.
Let $\Lambda$ be a dominant, integral
weight for $\Lhg$, then we can write  $\Lam= r \Lambda_0 + \lambda$
with $\lambda$ dominant, integral for $\Lg$.

Let $W$ be the $\Lg$-module $W := \overline{V}_{- \theta^{\vee}}(r
\Lam_0)$, where $\theta$ is the highest root of $\Lg$, we show that
W contains a unique one-dimensional submodule. Fix $w\neq 0$
a $\Lg$-invariant vector in $W$.
Let $V(\lambda)$ be the irreducible $\Lg$-module with highest weight $\lam$
and define the $\Lg$-module $V_{\lam,r}^{\infty}$ to be the direct limit of:
$$
V_{\lam,r}^\infty:\quad
V(\lambda) \hookrightarrow W \otimes V(\lambda)
\hookrightarrow W \otimes W \otimes V(\lambda)
\hookrightarrow W \otimes W \otimes W \otimes V(\lambda)
\hookrightarrow \ldots
$$
where the inclusions are always given by taking a
vector $u$ to its tensor product $u\mapsto w\otimes u$
with the fixed $\Lg$-invariant vector in $W$.

Recall the notation $\overline{V}(\Lambda)$ for $V(\Lam)$ viewed as
a $\Lg$-module. (The corresponding version
for a twisted Kac-Moody algebra can be found in section 4.)
\vskip 8pt\noindent {\bf Theorem 3.} {\it For any integral dominant
weight $\Lambda$ of $\Lhg$, $\Lam=r \Lambda_0 + \lambda$,
the $\Lg$-modules $V_{\lam,r}^{\infty}$ and $\overline{V}(\Lambda)$
are isomorphic.}
\vskip 8pt\noindent%
\begin{rem}\rm
The choice of $W$ is convenient because it avoids case by case
considerations. But, in fact, one could choose any other module
$W=V_{-\mu^{\vee}}( r \Lam_0)^{\otimes m}$, where $V_{-\mu^{\vee}}( r \Lam_0)$
is the Demazure module for a dominant, integral, nonzero coweight $\mu^{\vee}$ and $m$ is
such that $V_{-\mu^{\vee}}( r \Lam_0)^{\otimes m}$ contains a
one-dimensional submodule.
\end{rem}
\section{The affine Kac--Moody algebra}\label{affineKacMoody}
\subsection{Notations and basics}
In this section we fix the notation and the usual technical padding.
Let $\Lg$ be a simple complex Lie algebra. We fix a Cartan subalgebra
$\Lh$ in $\Lg$ and a Borel subalgebra $\Lb\supseteq \Lh$. Denote
$\Phi\subseteq \Lh^*$ the root system of $\Lg$, and, corresponding to the choice
of $\Lb$, let $\Phi^+$ be the set of positive roots
and let $\Delta=\{\al_1,\ldots,\al_n\}$ be the corresponding basis of $\Phi$.

For a root $\beta\in\Phi$ let $\beta^\vee\in\Lh$ be its coroot.
The basis of the dual root system (also called the coroot system)
$\Phi^\vee\subset \Lh$ is denoted $\Delta^\vee=\{\al_1^\vee,\ldots,\al_n^\vee\}$.

We denote throughout the paper by $\Theta=\sum_{i=1}^n a_i\al_i$ the
highest root of $\Phi$, by $\Theta^\times=\sum_{i=1}^n a_i^\times\al^\vee_i$
the highest root of $\Phi^\vee$
and by $\Theta^\vee=\sum_{i=1}^n a_i^\vee\al_i^\vee$
the coroot of $\Theta$. Note that $\Theta^\vee\not=\Theta^\times$ in general.
The Weyl group $W$ of $\Phi$ is generated
by the simple reflections $s_i=s_{\al_i}$ associated to the simple roots.

Let $P$ be the weight lattice of $\Phi$ and let $P^\vee$ be the
weight lattice of the dual root system $P^\vee$. Denote $P^+\subset P$
the subset of dominant weights and let $\bz[P]$ be the group algebra of $P$.
For a simple root $\al_i$ let $\om_i$ be the corresponding
fundamental weight, we use the same notation for simple coroots and
coweights. Recall that $\om_i$ is called {\it minuscule}
if $a^\times_i=1$, and the coweight $\om^\vee_i$ is called {\it minuscule}
if $a_i=1$.

Denote by $\Lh_\br\subset \Lh$
the real span of the coroots and let $\Lh_\br^* \subset \Lh^*$ be the real span of the
fundamental weights. We fix a $W$--invariant scalar product $(\cdot,\cdot)$
on $\Lh$ and normalize it such that the induced isomorphism
$\nu:\Lh_\br\longrightarrow \Lh_\br^*$
maps $\Theta^\vee$ to $\Theta$. With the notation as above it follows
for the weight lattice $P^\vee$ of the dual root system $\Phi^\vee$ that
$$
\nu(\al_i^\vee)=\frac{a_i}{a_i^\vee}\al_i\quad{\rm and}\quad
\nu(\om_i^\vee)=\frac{a_i}{a_i^\vee}\om_i,\quad\forall\,i=1,\ldots,n.
$$
Let $\Lhg$ be the affine Kac--Moody algebra corresponding
to the extended Dynkin diagram of $\Lg$ (see \cite{Kac}, Chapter 7):
$$
\Lhg=\Lg\otimes_\bc\bc[t,t^{-1}]\oplus \bc K\oplus \bc d
$$
Here $d$ denotes the derivation $d=t\frac{d}{dt}$ and $K$ is
the canonical central element. The Lie algebra $\Lg$ is
naturally a subalgebra of  $\Lhg$.
In the same way, $\Lh$ and $\Lb$ are
subalgebras of the Cartan subalgebra $\Lhh$ respectively the Borel
subalgebra $\Lhb$ of $\Lhg$:
\begin{equation}\label{hdecomp}
\Lhh=\Lh\oplus\bc K\oplus\bc d,\quad
\Lhb=\Lb\oplus\bc K\oplus\bc d\oplus\Lg\otimes_\bc t\bc[t]
\end{equation}
Denote by $\widehat\Phi$ the root system of $\Lhg$ and
let $\wPhi^+$ be the subset of positive roots. The positive non-divisible
imaginary root in $\wPhi^+$ is denoted $\delta$. The simple roots
are $\wDelta=\{\al_0\}\cup\Delta$ where $\al_0=\delta-\Theta$.
Let $\Lam_0,\ldots,\Lam_n$ be the corresponding fundamental weights,
then for $i=1,\ldots,n$ we have
\begin{equation}\label{Lamdecomp}
\Lam_i=\om_i+a_i^\vee\Lam_0.
\end{equation}
The decomposition of $\Lhh$ in (\ref{hdecomp}) has its corresponding
version for the dual space $\Lhh^*$:
\begin{equation}\label{hveedecomp}
\Lhh^* = \Lh^*\oplus\bc\Lam_0\oplus\bc \delta,
\end{equation}
here the elements of
$\Lh^*$ are extended trivially, $\langle\Lambda_0,\Lh\rangle
=\langle\Lambda_0,d\rangle=0$ and $\langle\Lambda_0,K\rangle=1$,
and $\langle\delta,\Lh\rangle=\langle\delta,K\rangle=0$ and
$\langle\delta,d\rangle=1$. Let $\wDelta^\vee=\{\al_0^\vee,\al^\vee_1,
\ldots,\al^\vee_n\}\subset \Lhh$ be the corresponding basis of the
coroot system, then $\al_0^\vee=K-\Theta^\vee$.

Set $\Lhh_\br^*=\br\delta+\sum_{i=0}^n\br\Lam_i$,
by (\ref{hveedecomp}) and (\ref{Lamdecomp})
we have  $\Lh_\br^* \subseteq \Lhh_\br^*$.
The affine Weyl group  $\waff$ is generated by the
reflections $s_0, s_1, . . . , s_n$, where again
$s_i=s_{\al_i}$ for a simple root. The cone
$\widehat{C} = \{ \Lam\in\Lhh_\br^* |\langle\Lam,\al_i^\vee\rangle\ge 0,
i= 0, . . . , n\}$ is the  fundamental Weyl chamber for $\Lhg$.

We put a $\widehat\ $ on (almost) everything
related to $\Lhg$. Let $\ph$ be the weight lattice of $\Lhg$, let
$\ph^+$ be the subset of dominant weights and let
$\bz[\ph]$ be the group algebra of $\ph$. Recall the following properties
of $\delta$ (see for example \cite{Kac}, Chapter 6):
\begin{equation}\label{Kaczitat}
\langle\delta,\al_i^\vee\rangle=0\,\forall\, i=0,\ldots,n\quad   w(\delta)=\delta
\,\forall\, w\in \waff,\quad
\langle\al_0,\al_i^\vee\rangle=-\langle\Theta,\al_i^\vee\rangle\,{\rm for\ }i\ge 1
\end{equation}
Put $a_0=a_0^\vee=1$ and let $A=(a_{i,j})_{0\le i,j\le n}$ be the (generalized)
Cartan matrix of $\Lhg$. We have a non--degenerate symmetric bilinear
form $(\cdot,\cdot)$ on $\Lhh$ defined by (\cite{Kac}, Chapter 6)
\begin{equation}\label{hatform}
\left\{
\begin{array}{ll}
(\al_i^\vee, \al_j^\vee )=\frac{a_j}{a_j^\vee}a_{i,j}&i,j=0,\ldots,n\\
(\al_i^\vee,d)=0&i=1,\ldots,n\\
(\al_0^\vee,d)=1&(d,d)=0.\\
\end{array}\right.
\end{equation}
The corresponding isomorphism $\nu:\Lhh\rightarrow\Lhh^*$ maps
$$
\nu(\al_i^\vee)=\frac{a_i}{a_i^\vee}\al_i,\quad
\nu(K)=\delta,\quad \nu(d)=\Lambda_0.
$$
Denote by $\Lg_{\rm sc}$ the subalgebra of $\Lhg$ generated by $\Lg$
and $\al_0^\vee=K-\Theta^\vee$, then $\Lh_{\rm sc}=\Lh\oplus\bc K$
is a Cartan subalgebra of $\Lg_{\rm sc}$. The inclusion $\Lh_\semc\rightarrow\Lhh$
induces an epimorphism $\Lhh^*\rightarrow \Lh_\semc^*$ with one dimensional
kernel. Now (\ref{Kaczitat}) implies that we have in fact an isomorphism
$$
\Lhh^*/\bc\delta\rightarrow \Lh_\semc^*
$$
and we set $\Lh_{\semc,\br}^*=\Lhh_\br^*/\br\delta$.
Since $\br\delta\subset \widehat{C}$, we use the same notation
$\widehat{C}$ for the image in $\Lh_{\semc,\br}^*$.
In the following we are mostly interested in characters of $\Lg$--modules
respectively $\Lg_\semc$--modules obtained by restricion from $\Lhg$--modules,
so we consider also the ring
$$
\bz[P_\semc]:=\bz[\ph]/I_\delta,
$$
where $I_\delta=(1-e^\delta)$ is the ideal in $\bz[\ph]$ generated by $(1-e^\delta)$.
\subsection{The extended affine Weyl group}
Since $\waff$ fixes $\delta$, the group can be defined as the
subgroup of $GL(\Lh_{\semc,\br}^*)$
generated by the induced reflections $s_0,\ldots,s_n$.

Let $M\subset \Lh^*_\br$
be the lattice $M=\nu(\bigoplus_{i=1}^n\bz\al_i^\vee)$. If $\Lg$ is simply
laced, then $M$ is the root lattice in $\Lh^*_\br$, otherwise $M$ is the lattice
in $\Lh^*_\br$ generated by the long roots.
An element $\Lam\in \Lh_{\semc,\br}^*$ can be uniquely decomposed
into $\Lam=\lam +b\Lam_0$ such that $\lam\in \Lh^*_\br$.
For an element $\mu\in M$ let $t_\mu\in GL(\Lh_{\semc,\br}^*)$ be the map
defined by
\begin{equation}\label{muaction}
\Lam =\lam +b\Lam_0\mapsto
t_\mu(\Lam)=\lam +b\Lam_0 + b\mu=\Lam+\langle \Lam, K\rangle\mu.
\end{equation}
Obviously we have $t_{\mu}\circ t_{\mu'}=t_{\mu+\mu'}$, denote $t_M$ the
abelian subgroup of $GL(\Lh_{\semc,\br}^*)$ consisting of the elements
$t_\mu$, $\mu\in M$. Then $\waff$ is the semidirect product
$\waff= W\semi t_M$.

The {\it extended affine Weyl group} $\wiff$ is the semidirect product
$\wiff= W\semi t_L$, where $L=\nu(\bigoplus_{i=1}^n\bz\om^\vee_i)$
is the image of the coweight lattice. The action of an element $t_\mu$,
$\mu\in L$, is defined as above in (\ref{muaction}).

Let $\Sigma$ be the subgroup of $\wiff$ stabilizing the dominant Weyl chamber
$\widehat{C}$:
$$
\Sigma=\{\sigma\in\wiff\mid \sigma(\widehat{C})=\widehat{C}\}.
$$
Then $\Sigma$ provides a complete system of coset representatives of
$\wiff/\waff$ and $\wiff=\Sigma\semi\waff$.
The elements $\sigma\in\Sigma$ are all of the form (one can verify
this easily or see \cite{Bourbaki})
$$
\sigma=\tau_i t_{-\nu(\omega_i^\vee)}= \tau _i t_{-\omega_i},
$$
where $\omega_i^\vee$ is a minuscule coweight and
$\tau_i=w_0 w_{0,i}$, where $w_0$ is the longest word in $W$
and $w_{0,i}$ is the longest word in $W_{\om_i}$, the stabilizer
of $\om_i$ in $W$.

We extend the length function $\ell:\waff\rightarrow\bn$ to a length
function $\ell:\wiff\rightarrow \bn$ by setting $\ell(\sigma w)=\ell(w)$
for $w\in \waff$ and $\sigma\in \Sigma$.%
\section{Demazure modules}
\subsection{Definitions}
For a dominant weight $\Lam\in \ph^+$ let $V(\Lam)$
be the (up to isomorphism) unique irreducible
$\Lhg$--highest weight module of highest weight $\Lam$.

Let $U(\Lhb)$ be the enveloping algebra of the Borel
subalgebra $\Lhb\subset\Lhg$. Given an element
$w\in\waff/W_\Lam$, fix a generator $v_{w(\Lam)}$
of the line $V(\Lam)_{w(\Lam)}=\bc  v_{w(\Lam)}$
of $\Lhh$--eigenvectors in $V(\Lam)$ of weight $w(\Lam)$.
\begin{dfn}\rm
The $U(\Lhb)$--submodule $V_w(\Lam)=U(\Lhb)\cdot v_{w(\Lam)}$
generated by $ v_{w(\Lam)}$ is called the {\it Demazure submodule
of $V(\Lam)$} associated to $w$.
\end{dfn}
To associate more generally to every element $\sigma w\in\wiff=\Sigma\semi\waff$
a Demazure module, recall that elements in $\Sigma$ correspond to automorphisms
of the Dynkin diagram of $\Lhg$, and thus define an associated automorphism
of $\Lhg$, also denoted $\sigma$. For a module $V$ of $\Lhg$
let $V^\sigma$ be the module with the twisted action $g\circ v=\sigma^{-1}(g)v$.
Then for the irreducible module of highest weight $\Lam\in \ph^+$ we get
$V(\Lam)^\sigma=V(\sigma(\Lam))$.

So for $\sigma w\in \wiff=\Sigma\semi\waff$ we set
\begin{equation}\label{demazuredefn}
V_{w\sigma }(\Lam)=V_{w}(\sigma(\Lam))
\quad{\rm respectively}\quad
V_{\sigma w}(\Lam)=V_{\sigma w\sigma^{-1}}(\sigma(\Lam)).
\end{equation}
Recall that for a simple root $\al$ the Demazure module $V_{w\sigma }(\Lam)$
is stable for the associated subalgebra ${\mathfrak{sl}}_2(\al)$ if
and only if $s_\al w\sigma\le w\sigma\mod W_\Lam$ in the (extended)
Bruhat order. In  particular, $V_{w\sigma }(\Lam)$ is a $\Lg$--module
if and only if $s_i w\sigma\le w\sigma\mod W^{\rm aff}_\Lam$ for all $i=1,\ldots,n$.

The example which will interest us are the Demazure modules
associated to the weight $r\Lam_0$ for $r\ge 1$, in this case $W^{\rm aff}_\Lam=W$,
so $\wiff/W=L$. The Demazure module  $V_{t_{\nu(\mu^\vee)}}(\Lam_0)$
is a $\Lg$--module if and only if $\mu^\vee$ is an anti-dominant coweight,
or, in other words, $\mu^\vee=-\lam^\vee$ for some dominant coweight.

To simplify the notation, we write in the following
\begin{equation}\label{nummereins}
V_{-\lam^\vee}(m\Lam_0)\quad{\rm for}\quad
V_{t_{-\nu(\lam^\vee)}}(m\Lam_0),
\end{equation}
and we write
\begin{equation}\label{nummerzwei}
\overline{V}_{-\lam^\vee}(m\Lam_0),
\end{equation}
for $V_{-\lam^\vee}(m\Lam_0)$ viewed as a $\Lg$--module. So we
view $\charc \overline{V}_{-\lam^\vee}(m\Lam_0)$ as an element in
$\bz[P]$ obtained from the $\Lhh$--character by projection.%
\subsection{Demazure operators} 
Let $\beta$ be a real root of the root system $\widehat\Phi$. We define
the {\it Demazure operator:}
$$
D_{\beta}:\bz[\ph]\rightarrow\bz[\ph],\quad
D_{\beta}(e^\lam)=\frac{e^\lam-e^{s_\beta(\lam)-\beta}}{1-e^{-\beta}}
$$
\begin{lem}\label{demazurezerleg}
\begin{enumerate}
\item For $\lambda,\mu\in\ph$ we have:
\begin{equation}\label{DemazureExpl}
D_{\beta}(e^\lam)=\left\{
\begin{array}{ll}
e^\lam+e^{\lam-\beta}+\dots+e^{s_\beta(\lam)}&
\hbox{\rm if\ }\langle\lam,\beta^\vee\rangle\ge 0\\
0&\hbox{\rm if\ }\langle\lam,\beta^\vee\rangle=-1\\
-e^{\lam+\beta}-e^{\lam+2\beta}-\dots-e^{s_\beta(\lam)-\beta}&
\hbox{\rm if\ }\langle\lam,\beta^\vee\rangle\le -2\\
\end{array}
\right.
\end{equation}
\item $D_{\beta}(e^{\lambda+\mu}) = e^{\lambda} D_{\beta}(e^{\mu}) +
e^{s_{\beta}(\mu)}D_{\beta}(e^{\lambda}) - e^{\lambda + s_{\beta}(\mu)}$ \\
\item Let $\chi \in \bz[\ph]$ be such that $s_{\beta}(\chi) = \chi $, then
$D_{\beta}(\chi) = \chi $.\\
\item Let $\chi \in \bz [\ph]$, then $D_{\beta}(\chi)$ is stable under $s_\beta$.
In particular, if $D_{\beta}(\chi) = \chi $, then
$s_{\beta}(\chi) = \chi $.\\
\item $D_{\beta}$ is idempotent, i.e.,
$D_{\beta}(D_{\beta}(e^{\mu})) = D_{\beta}(e^{\mu})$ for all $\mu$
\end{enumerate}
\end{lem}
\proof For 1., 3., 4., and 5. see \cite{Demazure}, (1.5)--(1.8).
The proof of  part 2. is a simple calculation.
\endpf
\vskip 5pt
Lemma~\ref{demazurezerleg} implies:
\begin{coro}
If $\langle\mu,{\beta}^\vee\rangle=0$, then
$D_{\beta}(e^{\lambda + \mu })=e^{\mu}D_{\beta}(e^{\lambda})$.
\end{coro}
The corollary is in fact a special case of the following more
general exchange rule, which follows easily from Lemma 1:
\begin{lem}\label{chimalcharformel}
Let $\chi,\eta\in\bz[\ph]$. If $D_{\beta}(\eta)=\eta$, then
$$D_{\beta}(\chi\cdot\eta)=\eta\cdot(D_{\beta}(\chi)).$$
\end{lem}
Since $D_{\alpha_i}(1-e^\delta)=(1-e^\delta)$ for all $i=0,\ldots,n$,
Lemma~\ref{chimalcharformel} shows that the ideal $I_\delta$ is
stable under all Demazure operators $D_{\beta}$. Thus we obtain
induced operators (we still use the same notation $D_\beta$)
$$
D_{\beta}:\bz[P_\semc]\longrightarrow \bz[P_\semc],\quad
e^\lam + I_\delta\mapsto D_{\beta}(e^\lam) + I_\delta,
$$
Recall further that $\langle\delta,{\beta}^\vee\rangle=0$ (see (\ref{Kaczitat})),
so it makes sense to define on $\bz[P_\semc]$ the function
$e^\lam\mapsto \langle\lam,{\beta}^\vee\rangle$.
\begin{lem}\label{alphanullundtheta}
If $\lam\in \ph\cap\Lh^*$, then $D_{\al_0}(e^\lam)=D_{-\Theta}(e^\lam)$
in $\bz[P_\semc]$.
\end{lem}
\proof
Since $\lam\in\Lh^*$ we have $\langle\lam,\al_0^\vee\rangle
=\langle\lam,c-\Theta^\vee\rangle=-\langle\lam,\Theta^\vee\rangle$.
Further, $\al_0=\delta-\Theta$, so equation~(\ref{DemazureExpl})
can be read in $\bz[P_\semc]$ as
\begin{equation}\label{DemazureExpl3}
\begin{array}{rcl}
D_{\al_0}(e^\lam)&=&\left\{
\begin{array}{ll}
e^\lam+e^{\lam+\Theta}+\dots+e^{\lam+n\Theta}&
\hbox{\rm if\ }n=\langle\lam,\al_0^\vee\rangle=\langle\lam,-\Theta^\vee\rangle\ge 0\\
0&\hbox{\rm if\ }\langle\lam,\al_0^\vee\rangle=\langle\lam,-\Theta^\vee\rangle=-1\\
-e^{\lam-\Theta}-\dots-e^{\lam-(\vert n\vert-1)\Theta}&
\hbox{\rm if\ }n=\langle\lam,\al_0^\vee\rangle=\langle\lam,-\Theta^\vee\rangle\le -2\\
\end{array}
\right.\\
\\
&=&D_{-\Theta}(e^\lam)\hfill\bullet\\
\end{array}
\end{equation}
\subsection{Demazure character formula}
We want to extend the notion of a Demazure operator also to
elements of $\Sigma$. We define for $\sigma\in \Sigma$:
$$
D_\sigma: \bz[\ph]\rightarrow \bz[\ph],\quad
D_\sigma(e^\Lam)=e^{\sigma(\Lam)}.
$$
Since $\sigma(\delta)=\delta$, we get an induced operator $D_\sigma$
on $\bz[P_\semc]$.
\begin{lem}\label{sigmaoperator}
$D_\sigma D_\beta=D_{\sigma(\beta)}D_\sigma$.
\end{lem}
\proof
Let $\Lam\in \ph$, then $\langle\Lam,\beta^\vee\rangle=
\langle \sigma(\Lam),\sigma(\beta^\vee)\rangle$, which implies
the claim by equation ~(\ref{DemazureExpl}).
\endpf
In the following we denote by $D_i$, $i=0,\ldots,n$ the Demazure operator
$D_{\al_i}$ corresponding to the simple root $\al_i$.
Recall that for any reduced decomposition $w=s_{i_1}\cdots s_{i_r}$
of $w\in\waff$ the operator $D_w=D_{i_1}\cdots D_{i_r}$
is independent of the choice of the decomposition (see \cite{Kumar},
Corollary 8.2.10).

We associate an operator to any element $w\sigma \in\wiff$
by setting
$$
\begin{array}{rcl}
D_{w\sigma}:\bz[P_\semc]&\rightarrow&\bz[P_\semc]\\
e^\Lam&\mapsto&D_w(e^{\sigma(\Lam)})\\
\end{array}
$$
By Lemma~\ref{sigmaoperator} we have for $\sigma w\in\wiff=\Sigma\semi\waff$:
$$
\begin{array}{rcl}
D_{\sigma w}:\bz[P_\semc]&\rightarrow&\bz[P_\semc]\\
e^\Lam&\mapsto&\sigma\Big(D_w(e^\Lam)\Big)=
D_{\sigma w\sigma^{-1}}(e^{\sigma(\Lam)})\\
\end{array}
$$
Let $w\sigma\in\wiff$ and let $\Lam\in \ph^+$ be a dominant weight.
\begin{thm}[\cite{Kumar} Chapter VIII, \cite{Kumar2,Mathieu}]\label{demazurcharacterformula}
$$
\charc V_w(\sigma(\Lam))=D_{w\sigma}(e^\Lam).
$$
\end{thm}
Let $\lam^\vee$ be a dominant coweight. Associated to
$t_{-\nu(\lam)^\vee}\in \wiff$ we have a Demazure operator
$D_{t_{-\nu(\lam^\vee)}}$, we write for simplicity just
$D_{-\lam^\vee}$.
\begin{lem}\label{demazuradd}

Let $\lam_1^\vee,\lam_2^\vee$ be two dominant coweights, and set $\lam^\vee=
\lam_1^\vee+\lam_2^\vee$. Then
$$
D_{-\lam_1^\vee}D_{-\lam_2^\vee}=D_{-\lam^\vee}
$$
\end{lem}
\begin{lem}\label{sechs}
Let $V$ be a finite dimensional
$\Lg_\semc$--module such that $\charc V\in \bz[P]$, then
\begin{equation}\label{quotientformel}
D_i(\charc V)=\charc V\quad\forall i=0,\ldots,n;
{\rm\ and\ }D_\sigma(\charc V)=\charc V.
\end{equation}
\end{lem}
\proof
The character of a finite dimensional $\Lg$--module is stable under the Weyl group
$W$ and hence stable under $D_i$ for all $i=1,\ldots,n$ by Lemma~\ref{demazurezerleg}.
It remains to consider the case  $i=0$. Now all weights lie in $\Lh^*$,
so by Lemma~\ref{alphanullundtheta} we have:
$$
D_0(\charc V)=D_{-\Theta}(\charc V)=\charc V
$$
where the right hand side is again a consequence of
Lemma~\ref{demazurezerleg}.

Now $\sigma=y t_{-\nu(\om_j^\vee)}$ for some minuscule fundamental
coweight $\om_j^\vee$ and some $y\in W$. Since $
t_{-\nu(\om_j^\vee)}$ operates trivially on $\bz[P]$ and
$D_y(\charc V)=\charc V$, the claim follows.
\hfill $\bullet$%
\section{The proofs}
\subsection{Proof of Theorem 1}
Let $\lam^\vee$ be a dominant coweight and suppose we are given a
decomposition
$$
\lam^\vee=\lam_1^\vee+\lam_2^\vee+\ldots+\lam_r^\vee
$$
of $\lam^\vee$ as a sum of dominant coweights. For the notation
see (\ref{nummereins}) and (\ref{nummerzwei}).%
\vskip 5pt \noindent%
{\bf Theorem 1} \label{Theoremeins} {\it As
$\Lg$--representations, the modules }
$$
\overline{V}_{{-\lam^\vee}}(m\Lam_0)\quad{\rm and}\quad
\overline{V}_{{-\lam_1^\vee}}(m\Lam_0)\otimes
\overline{V}_{{-\lam_2^\vee}}(m\Lam_0) \otimes\cdots\otimes
\overline{V}_{{-\lam_r^\vee}}(m\Lam_0)
$$
{\it are isomorphic.}\\ \\%
More precisely, we will show that, on the level of characters of
$\Lg_\semc$--modules:%
\vskip 8pt\noindent%
{\bf Theorem 1'.}
$$
\charc {V}_{{-\lam^\vee}}(m\Lam_0)= e^{m\Lam_0} \charc
\overline{V}_{{-\lam_1^\vee}}(m\Lam_0) \charc
\overline{V}_{{-\lam_2^\vee}}(m\Lam_0) \cdots \charc
\overline{V}_{{-\lam_r^\vee}}(m\Lam_0).
$$
Theorem 1' obviously implies Theorem 1, so
it suffices to prove Theorem 1'.\\
A first step is the following lemma:
\begin{lem}\label{firststep}
Let $\chi\in \bz[P_\semc]$ be a character of the form
$e^{m\Lam_0}\charc \overline{V}$, where $\overline{V}$ is a finite
dimensional $\Lg$--module. Suppose $\lam^\vee\in P^\vee$ is a
dominant coweight and let $t_{-\nu(\lam^\vee)}= s_{i_1}\dots
s_{i_t}\sigma$ be a reduced decomposition in $\wiff$. Then
$$
D_{{i_1}}\dots D_{i_t} D_\sigma(e^{m\Lam_0}\charc \overline{V})=
D_{{i_1}}\dots D_{i_t} D_\sigma(e^{m\Lam_0})\charc \overline{V}.
$$
\end{lem}
\proof The lemma is proven exactly in the same way as
Lemma~\ref{chimalcharformel}, only using now in addition
Lemma~\ref{sechs} for the operators $D_0$ and $D_\sigma$. \hfill
$\bullet$ \vskip 5pt\noindent {\it Proof of
Theorem~\ref{Theoremeins}'.} The proof is by induction on $r$.
Suppose $r=1$ and $\lam^\vee=w\sigma$ where $\sigma\in \Sigma$ and
$w\in\waff$. The character of $V_{-\lam^\vee}(m\Lam_0)$ is the
character of the Demazure submodule
$V_w(\sigma(m\Lam_0))=V_w(m\Lam_0 + m\om_i^*)$ for some
appropriate minuscule fundamental weight of $\Lg$. So all
$\Lg_\semc$--weights occuring in the module are of the form
$m\Lam_0+m\om_i^*+$ a sum of roots in $\Phi$ (possibly positive
and negative, see Lemma 3), and hence the character is of the
desired form $e^{m\Lam_0}\charc
\overline{V}_{-\lam^\vee}(m\Lam_0)$.

Suppose now $r\ge 2$ and the claim holds already for $r-1$. By
the definition in equation~(\ref{demazuredefn}) we have for
$t_{\nu(-\lam^\vee)}=w\sigma\in\wiff$:
$$
\charc V_{-\lam^\vee}(m\Lam_0)=\charc V_{w}(m\sigma(\Lam_0)),
$$
by the Demazure character formula
(Theorem~\ref{demazurcharacterformula}) the latter is equal to
$D_{-\lam^\vee}(e^{m\Lam_0})$, so
$$
\charc V_{-\lam^\vee}(m\Lam_0)=D_{-\lam^\vee}(e^{m\Lam_0}),
$$
by Lemma~\ref{demazuradd} the right hand side can be rewritten as
$$
\charc V_{-\lam^\vee}(m\Lam_0)=D_{-\lam_1^\vee}
\Big(D_{-\lam_2^\vee}\cdots D_{-\lam_r^\vee}(e^{m\Lam_0})\Big),
$$
by induction the right hand side can be reformulated as
$$
\charc V_{-\lam^\vee}(m\Lam_0)=D_{-\lam_1^\vee}\Big( e^{m\Lam_0}
\charc \overline{V}_{{-\lam_2^\vee}}(m\Lam_0) \cdots \charc
\overline{V}_{{-\lam_r^\vee}}(m\Lam_0)\Big),
$$
by Lemma~\ref{firststep} this is equivalent to
$$
\charc V_{-\lam^\vee}(m\Lam_0)=
\Big(D_{-\lam_1^\vee}(e^{m\Lam_0})\Big) \charc
\overline{V}_{{-\lam_2^\vee}}(m\Lam_0) \cdots \charc
\overline{V}_{{-\lam_r^\vee}}(m\Lam_0).
$$
Now the arguments for the proof of the case $r=1$ show that this implies
$$
\charc V_{-\lam^\vee}(m\Lam_0)=
e^{m\Lam_0}\charc \overline{V}_{{-\lam_1^\vee}}(m\Lam_0)
\charc \overline{V}_{{-\lam_2^\vee}}(m\Lam_0)
\cdots
\charc \overline{V}_{{-\lam_r^\vee}}(m\Lam_0),
$$
which finishes the proof.
\endpf
\subsection{Proof of Theorem 1 A} The proof is similar to the proof
above, so we give just a short sketch. As above, we have
$$
\charc V_{-\lam^\vee}(m\Lam_0+r\Lam_i)
=D_{-\lam^\vee} \Big(e^{m\Lam_0+r\Lam_i}\Big)
=D_{-\lam_2^\vee}D_{-\lam_3^\vee}\cdots D_{-\lam_r^\vee}
D_{-\om_i^\vee}\Big(e^{m\Lam_0+r\Lam_i}\Big).
$$
Now $t_{-\nu(\om_i)}=t_{-\om_i}=\tau_i\sigma_i$. Here
$\tau_i= w_{0,i}w_0$, where $w_0$ is the longest element in $W$
and $w_{0,i}$ is the longest word in the stabilizer $W_{\om_i}$ of $\om_i$,
and $\sigma_i$ is a diagram automorphism. Note that
$\sigma_i(\Lam_i)=\tau_i^{-1}t_{-\om_i}(\Lam_0+\om_i)
=\tau_i^{-1}(\Lam_0)=\Lam_0$ and
$$
\sigma_i(\Lam_0)=\tau_i^{-1}t_{-\om_i}(\Lam_0)
=\tau_i^{-1}(\Lam_0-\om_i)=\Lam_0+\tau_i^{-1}(-\om_i)=
\Lam_0+w_{0}w_{0,i}(-\om_i)=\Lam_0+\om_i^*,
$$
where $\omega^*$ denotes the highest weight of the irreducible
$\Lg$-representation $V(\om_i)^*$. Note that $\Lam_0+\om_i^*$
is again a fundamental weight (for the Kac--Moody algebra $\Lhg$),
and recall that
$$
\tau_i=w_{0,i}w_0
=w_{0}(w_0^{-1}w_{0,i}w_0)=w_0w_{0,i}^*,
$$
where $w_{0,i}^*$ is the longest word in the stabilizer
$W_{\om_i^*}$ of $\om_i^*$. So
$$
\begin{array}{rcl}
D_{-\om_i^\vee}\Big( e^{m\Lam_0+r\Lam_i} \Big)&=&
D_{\tau_i} D_{\sigma_i} \Big( e^{m\Lam_0+r\Lam_i} \Big)\\
&=&D_{\tau_i}\Big( e^{m\Lam_0+m\om_i^* +r\Lam_0} \Big)\\
&=& e^{(m+r)\Lam_0} D_{\tau_i}\Big( e^{m\om_i^*}\Big)\\
&=& e^{(m+r)\Lam_0} D_{w_0w_{0,i}^*}\Big( e^{m\om_i^*}\Big)\\
&=& e^{(m+r)\Lam_0} \charc{V}(m\om_i^*).\\
\end{array}
$$
Now the same induction procedure as above applies to finish the proof.
\endpf
\subsection{Proof of Theorem 2}
The proof is divided into several
case by case considerations. Suppose first that $\om^\vee$ is a
minuscule coweight. In this case (for $m=1$ this has already been
proved in \cite{Magyar}) $t_{-\om_i}=w_{i,0}w_0\sigma_i$ and hence
$$
\charc \overline{V}_{-\om^\vee}(m\Lam_0)=D_{-\om^\vee}e^{m\Lam_0}
=D_{w_{i,0}w_0}e^{m\Lam_0+m\om^*}
=D_{w_0w^*_{i,0}}e^{m\Lam_0+m\om^*}
=e^{m\Lam_0}\charc V(m\om)^*.
$$
In particular, this finishes the proof for the Lie algebras of type ${\tt A}_n$.
For the next few cases we need the following:
\begin{lem}\label{hilf8}
Let $w_0$ be the longest element in the Weyl group of $\Lg$, let $z$ be an arbitrary element
of $\mbox{Stab}_{W}(\Theta)$, where $\Theta$ is the highest root of $\Lg$, let $r \in \bn$. Then\\
$$
\overline{V}_{w_0 z s_0}(r \Lam_0) \simeq \bigoplus\limits_{m = 0}^{r} V(m \Theta)
$$
as $\Lg$-representations.
\end{lem}
\proof
$$
\begin{array}{rcl}
D_{(w_0 z) s_0}(e^{r \Lam_0}) &= & D_{w_0 z}D_{-\Theta}(e^{r \Lam_0}) \\
& =& D_{w_0 z}( e^{r \Lam_0} + e^{r \Lam_0 + \Theta} + \ldots + e^{ r \Lam_0 + r \Theta})\\
& =& D_{w_0}( e^{r \Lam_0} + e^{r \Lam_0 + \Theta} + \ldots + e^{ r \Lam_0 + r \Theta})\\
& =& e^{r \Lam_0}(D_{w_0}( e^{0} + D_{w_0}( e^{\Theta}) + \ldots +  D_{w_0}( e^{r \Theta})))\\
\end{array}
$$
which finishes the proof.
\endpf
In the cases $E_6, E_7, E_8, F_4, G_2$ the highest root $\Theta$
is also a fundamental weight, say $\om_i$. Let $p_i:= \frac{a_i
}{a^{\vee }_i}$. Then $\nu(\om_i^{\vee}) = \frac{a_i}{a^{\vee}_i}
\om_i = p_i\Theta$. In fact, for the adjoint representations considered
here one sees that $ p_i=1$ in all cases. Since $t_{-\om_i}=s_\theta s_0$,
it follows by Lemma~\ref{hilf8}:
$$
\overline{V}_{-\om_i^{\vee}}(r \Lam_0) \simeq \bigoplus\limits_{m
= 0}^{r} V(m \om_i)
$$
Next we consider the types ${\tt B}_n$ and ${\tt D}_n$ with the Bourbaki indexing of the simple
roots, i.e., we consider the root system as embedded in $\br^n$ with the canonical basis
$\{\eps_1,\ldots,\eps_n\}$ and the standard scalar product.
The basis of the root system is given by the simple roots $\al_i=\eps_i-\eps_{i+1}$,
$i=1,\ldots,n-1$ and $\al_n=\eps_n$ (type ${\tt B}_n$, $n\ge 3$) respectively
$\al_n=\eps_{n-1}+\eps_{n}$ (type ${\tt D}_n$, $n\ge 4$), the highest root is
$\eps_1+\eps_2$ in both cases. We have
$$
\begin{array}{rcl}
t_{-\om_2}&=&s_{\eps_1+\eps_2}s_0\\
&=&(s_2\cdots s_n\cdots s_2)s_1
(s_2\cdots s_n\cdots s_2)s_0\\
\end{array}
$$
In the following we consider only the non-minuscule fundamental coweights.
We get for $2i\le n$ (case ${\tt B}_n$) respectively $2i\le n-2$ (case ${\tt D}_n$):
\begin{equation}\label{last}
\begin{array}{rcl}
t_{-\nu(\om_{2i}^\vee)}\hskip -8pt
&=&\hskip -8ptt_{-\eps_1-\eps_2}t_{-\eps_3-\eps_4}\cdots t_{-\eps_{2i-1}-\eps_{2i}}\\
&=&\hskip -8ptt_{-\om_2}\Big((s_2s_1s_3s_2)t_{-\om_2}(s_2s_1s_3s_2)\Big)\cdots
\Big((s_{2i-2}\cdots s_2s_1)\\
&&(s_{2i-1}\cdots s_3s_2)t_{-\om_2}
(s_2s_3\cdots s_{2i-1})(s_1s_2\cdots s_{2i-2})\Big)\\
&=&\hskip -8pt s_{\eps_1+\eps_2}s_0 s_{\eps_3+\eps_4}\Big((s_2s_1s_3s_2)s_0(s_2s_1s_3s_2)\Big)
s_{\eps_5+\eps_6}\cdots s_{\eps_{2i-1}+\eps_{2i}}\\
&&
\Big((s_{2i-2}\cdots s_2s_1)(s_{2i-1}\cdots s_3s_2)s_0
(s_2s_3\cdots s_{2i-1})(s_1s_2\cdots s_{2i-2})\Big)\\
&=&\hskip -8pt\Big[s_{\eps_1+\eps_2}s_{\eps_3+\eps_4}\cdots s_{\eps_{2i-1}+{\eps_{2i}}}\Big]
\Big[s_0 \Big((s_2s_1s_3s_2)s_0(s_2s_1s_3s_2)\Big)\cdots\\
&&\hskip -8pt \Big((s_{2i-2}\cdots s_2s_1)(s_{2i-1}\cdots s_3s_2)s_0
(s_2s_3\cdots s_{2i-1})(s_1s_2\cdots s_{2i-2})\Big)\Big]\\
\end{array}
\end{equation}
We see that we can write the word as a product $w_1w_2$ of two words,
the first being an element of the Weyl group $W$ and the second being a word
in the subgroup of $\waff$ generated by the simple reflections
$s_0,s_1,\ldots,s_{2i-1}$, this is (in the ${\tt B}_n$ as well as in the ${\tt D}_n$
case) a group of type ${\tt D}_{2i}$.

Since we look for a character of a $\Lg$--module, we know the character is stable
under the operators $D_i$, $1\le i\le n$. So to determine
the character of $V_{-\om_{2i}^\vee}(m\Lam_0)$, it suffices to
get a reduced decomposition of the word $w_2$ above modulo the right and left action
of $W$, the character of $V_{-\om_{2i}^\vee}(m\Lam_0)$ can be reconstructed
by applying the Demazure operators $D_i$, $1\le i\le n$.

The strategy is the following. We show that the decomposition above of $w_2$ is
a reduced decompostion. Further, we show that $\tau=s_1 s_3\ldots s_{2i-1} w_2$
is the longest word of the Weyl group of the subgroup of $\waff$ of type ${\tt D}_{2i}$.

Before we give a more detailed account on how to prove this, let us show
how this solves the problem. Let $\Ld\subset \Lhg$ be the semisimple
Lie algebra of type ${\tt D}_{2i}$ associated to the simple roots
$\al_0,\ldots,\al_{2i-1}$, then $V_\tau(m\Lam_0)$ is
an irreducible $\Ld$--module. More precisely,
it is the irreducible $m$--th spin representation (associated to the node
of $\al_0$). Let $\Ld'$ be the semisimple subalgebra of $\Ld$ corresponding
to the simple roots $\al_1,\ldots,\al_{2i-1}$, then $\Ld'$ is also
the semisimple part of a Levi subalgebra of $\Lg$.
Since $V_\tau(m\Lam_0)$ is a $\Lhb$--module, it is hence a $\Lb$ and a $\Ld'$--module.
By the Borel-Weil-Bott theorem we know that the induced
$\Lg$--module (which is the module $V_{-\om_{2i}^\vee}(m\Lam_0)$)
has the same direct sum decomposition as
$V_\tau(m\Lam_0)$ has as $\Ld'$--module. Since the latter has
been already given in \cite{Littelmann}, this finishes the proof.

We come now back to the proof of the first claim.
We make  the calculations in the following modulo $\delta$,
so the set of positive roots for the type ${\tt D}_{2i}$--subdiagram
(modulo $\delta$) is the set
$$
\{\eps_s-\eps_t\mid 1\le s<t\le 2i \}\cup\{-\eps_s-\eps_t\mid 1\le s<t\le 2i\}.
$$
In these terms the decomposition of $w_2$ as the second part in the
square brackets in (\ref{last}), reads as
$$
w_2=s_{-\eps_1-\eps_2}s_{-\eps_3-\eps_4}\cdots s_{-\eps_{2i-1}-\eps_{2i}}.
$$
and all positive roots above are sent to negative by $w_2$ roots except
$\al_1,\al_3,\ldots,\al_{2i-1}$. This implies that the decomposition above
of length $4i^2-3i$ is reduced, and
\begin{equation}\label{dword}
\begin{array}{rl}
\tau=s_1 s_3\ldots s_{2i-1}s_0 \Big((s_2s_1s_3s_2)s_0(s_2s_1s_3s_2)\Big)
\hskip -9pt&\cdots\hskip -2pt\Big((s_{2i-2}\cdots s_2s_1)(s_{2i-1}\cdots s_3s_2)\\
&\quad \quad s_0 (s_2s_3\cdots s_{2i-1})(s_1s_2\cdots s_{2i-2})\Big)\\
\end{array}
\end{equation}
is a reduced decomposition of the longest word of the Weyl group of the subgroup
of type ${\tt D}_{2i}$. This shows that $\tau$ a subword of $t_{-\om_{2i}}$,
and $V_{-\om_{2i}^\vee}(m\Lam_0)$ is the $\Lg$--module
generated by the $\Lb$--$\Ld'$--submodule $V_\tau(m\Lam_0)$.

It has been already pointed out above that the decomposition
of $V_{-\om_{2i}^\vee}(m\Lam_0)$ as $\Lg$--module is completely
determined by the $\Lh$--module structure of $V_\tau(m\Lam_0)$
and the decomposition of $V_\tau(m\Lam_0)$ as $\Ld'$--module.
So it remains to describe the decompositon of the $m$-th spin--representation
$V_\tau(m\Lam_0)$ with respect to the subalgebra $\Ld'$, and to describe
the highest weights as weights for the Cartan subalgebra $\Lh$.

The decomposition of the $m$-th spin--representation
$V_\tau(m\Lam_0)$ with respect to the subalgebra $\Ld'$
can be found in \cite{Littelmann} (see section 1.4). The description
of the possible highest weights occuring (\cite{Littelmann}, Proposition 3.2)
in the decomposition implies for $2i<n$ (case ${\tt B}_n$)
respectively $2i\le n-2$ (case ${\tt D}_n$):
$$
\charc\overline{V}_{-\om_{2i}^\vee}(m\Lam_0)=\sum_{a_1+\ldots+a_{i}=m}
\charc V(a_1\om_2+\ldots+a_{i}\om_{2i}),
$$
and for $2i=n$ in the case ${\tt B}_n$:
$$
\charc\overline{V}_{-\om_{n}^\vee}(m\Lam_0)=\sum_{a_1+\ldots+a_{n/2}=m}
\charc V(a_1\om_2+\ldots+a_{(n-2)/2}\om_{n-2}+2a_n\om_n).
$$
The calculation for the odd case is similar. We assume $2i+1\le n$ in the case ${\tt B}_n$
and $2i+1\le n-2$ in the case ${\tt D}_n$
$$
\begin{array}{rcl}
t_{-\om_{2i+1}^\vee}&=&
t_{-\eps_1-\eps_2}\cdots t_{-\eps_{2i-1}-\eps_{2i}}t_{-\eps_{2i+1}}\\
&=&t_{-\eps_1-\eps_2}(s_2s_1s_3s_2t_{-\eps_1-\eps_2}s_2s_3s_1s_2)\cdots\\
&&\ (s_{2i-2}\ldots s_1s_{2i-1}\ldots s_2t_{-\eps_1-\eps_2}s_2\ldots s_{2i-1} s_1\ldots s_{2i-2})\\
&&\quad (s_{2i}\ldots s_1 t_{-\eps_1} s_1\ldots s_{2i})\\
&=&s_{\eps_1+\eps_2}s_0s_{\eps_3+\eps_4}(s_2s_1s_3s_2 s_0s_2s_3s_1s_2)\cdots\\
&&\ s_{\eps_{2i-1}+\eps_{2i}}(s_{2i-2}\ldots s_1s_{2i-1}\ldots s_2s_0s_2\ldots s_{2i-1} s_1\ldots s_{2i-2})\\
&&\quad s_{\eps_{2i+1}}(s_{2i}\ldots s_1 \sigma_1 s_1\ldots s_{2i})\\
&=&[s_{\eps_1+\eps_2}s_{\eps_3+\eps_4}\cdots s_{\eps_{2i-1}+\eps_{2i}}s_{\eps_{2i+1}}]
[s_0(s_2s_1s_3s_2 s_0s_2s_3s_1s_2)\cdots\\
&&\ (s_{2i-2}\ldots s_1s_{2i-1}\ldots s_2s_0s_2\ldots s_{2i-1} s_1\ldots s_{2i-2})
(s_{2i}\ldots s_1 s_0s_2\ldots s_{2i})\sigma_1]\\
\end{array}
$$
It follows as above that the second part of the word is reduced.  In fact,
after multiplying the word with $s_1s_3\ldots s_{2i-1}$, we obtain a reduced decomposition
of the longest word
$$
\begin{array}{rl}
\tau=s_1s_3\ldots s_{2i-1}s_0(s_2s_1s_3s_2 s_0s_2s_3s_1s_2)\cdots\hskip -2pt&\hskip -8pt
(s_{2i-2}\ldots s_1s_{2i-1}\ldots s_2s_0s_2\ldots s_{2i-1} s_1\ldots s_{2i-2})\\
&(s_{2i}\ldots s_1 s_0s_2\ldots s_{2i})\\
\end{array}
$$
in the Weyl group of the semisimple Lie algebra $\Ld\subset \Lhg$ of
type ${\tt D}_{2i+1}$
associated to the simple roots $\al_0,\ldots,\al_{2i}$. The Demazure module
$V_{\tau\sigma_1}(m\Lam_0)$ is an irreducible $\Ld$--module, it is the $m$--th
spin representation, associated to the node corresponding to $\al_1$. Consider
the decomposition of $V_{\tau\sigma_1}(m\Lam_0)$ as an $\Lh$-- and a
$\Ld'$--module, where $\Ld'\subset\Ld$ is the semisimple Lie subalgebra
associated to the simple roots $\al_1,\ldots,\al_{2i}$. By \cite{Littelmann},
we get as $\Ld'$--$\Lh$--module the decomposition
($2i+1<n$ in the ${\tt B}_n$ case):
$$
\overline{V}_{\tau\sigma_1}(m\Lam_0)=
\overline{V}_{\tau}(m\Lam_1)=
\bigoplus_{a_1+\ldots+ a_i= m}V(a_1\om_1+a_1\om_3+\ldots+a_i \om_{2i+1})
$$
and, again by the Borel--Weil--Bott theorem,  the same decomposition
holds for the Demazure module $V_{-\om_{2i+1}^\vee}(m\Lam_0)$ as $\Lg$--module.
The case $n=2i+1$ is treated similarly.

Next we consider the Lie algebra of type ${\tt C}_n$. We have for
$j=1,\ldots,n-1$ ($\omega_n^\vee$ is minuscule)
$$
t_{-\om_j^\vee}=t_{-2\om_j}=t_{-2\eps_1}t_{-2\eps_2}\cdots t_{-2\eps_j}=
t_{-2\eps_1}(s_1t_{-2\eps_1}s_1)\cdots (s_{j-1}\cdots s_1t_{-2\eps_1}s_1\cdots s_{j-1}).
$$
Replacing $t_{-2\eps_1}$ by $s_{2\eps_1}s_0$ we get
$$
\begin{array}{rcl}
t_{-\om_j^\vee}
&=&s_{2\eps_1}s_0(s_1s_{2\eps_1}s_0s_1)(s_2s_1s_{2\eps_1}s_0s_1s_2)
\cdots (s_{j-1}\cdots s_1s_{2\eps_1}s_0s_1\cdots s_{j-1})\\
&=&s_{2\eps_1}s_0s_{2\eps_2}(s_1s_0s_1)s_{2\eps_3}(s_2s_1s_0s_1s_2)
\cdots s_{2\eps_{j}}(s_{j-1}\cdots s_1s_0s_1\cdots s_{j-1})\\
&=&[s_{2\eps_1}s_{2\eps_2}s_{2\eps_3}
\cdots s_{2\eps_{j}}][
s_0(s_1s_0s_1)(s_2s_1s_0s_1s_2)\cdots
(s_{j-1}\cdots s_1s_0s_1\cdots s_{j-1})]\\
\end{array}
$$
We proceed now with the same strategy as before. For the moment we omit
the reflections $s_{2\eps_1}s_{2\eps_2}s_{2\eps_3}\cdots s_{2\eps_{j}}$.
The second part, the word
$$
\tau=s_0(s_1s_0s_1)(s_2s_1s_0s_1s_2)\cdots
(s_{j-1}\cdots s_1s_0s_1\cdots s_{j-1}),
$$
is a reduced decomposition of the longest
word of the semisimple subalgebra $\Ld\subset \Lhg$ of type ${\tt C}_{j}$
associated to the simple roots $\al_0,\ldots,\al_{j-1}$.

The Demazure module $V_\tau(m\Lam_0)$ is, as $\Ld$--module,
irreducible. Let $\Ld'\subset \Ld$ be the semisimple Lie algebra
associated to the simple roots $\al_1,\ldots,\al_{j-1}$,
it follows again from \cite{Littelmann} that the restriction of
$V_\tau(m\Lam_0)$ decomposes as $\Ld'$-- and $\Lh$--module
$$
\overline{V}_\tau(m\Lam_0)\simeq \bigoplus_{a_1+\ldots+a_j\le m}
\overline{V}(2a_1\om_1+\ldots +2a_j\om_j),
$$
which, as above, implies the corresponding decomposition as
$\Lg$-module.%

Fot $\Lg$ of type ${\tt F}_4$ and $\om_4^\vee$ we use the same strategy as above.
Using the same notation as in~\cite{Bourbaki}, one sees  $2\om_4=2\eps_1=(\eps_1+\eps_2)
+(\eps_1-\eps_2)=\Theta+s_1s_2s_3s_2s_1(\Theta)$, so
$$
\begin{array}{rcl}
t_{-\nu(\om_4^\vee)}=t_{-2\om_4}=t_{-\eps_1-\eps_2}t_{-\eps_1+\eps_2}
&=&(s_\Theta s_0)(s_1s_2s_3s_2s_1 s_\Theta s_1s_2s_3s_2s_1)\\
&=&(s_\Theta s_{\eps_1-\eps_2})(s_0 s_1s_2s_3s_2s_1s_0 s_1s_2s_3s_2s_1)\\
\end{array}
$$
Again we decompose the translation into a product of two words $w\tau$ such
that $w\in W$ and $\tau$ is a subword of a reduced decomposition of the
longest word of the Weyl group of type ${\tt B_4}$ corresponding to the
roots $\{\al_0,\al_1,\al_2,\al_3\}$. For the corresponding Levi subalgebra of $\Lhg$
we have $V_\tau(m\Lam_0)$ is the Cartan component in the $m$-th symmetric
power of the action of the orthogonal Lie algebra on $\bc^n$. Now by looking
at the decomposition of this space  with respect to the Levi subalgebra
of $\Lg$ corresponding to the simple roots $\{\al_1,\al_2,\al_3\}$ (using
again the tables in \cite{Littelmann}), we obtain the desired formula.
\subsection{Proof of Theorem 3}
We will need the following simple:
\begin{lem}\label{cogewicht+w0} Let $\lam^{\vee}$ be a dominant, integral
coweight of $\Lg$, let $w_0$ be the longest element of the
Weyl group of $\Lg$, then:
$$l(t_{-\lam^{\vee}}w_0)=l(t_{-\lam^{\vee}})+l(w_0)$$
\end{lem}
\noindent %
So reduced decompositions of $t_{-\lam^{\vee}}$ and $w_0$ give a
reduced decomposition of $t_{-\lam^{\vee}} w_0$.\\ \\%
\begin{lem}
Let $W$ be the $\Lg$-module $W := \overline{V}_{-\theta^{\vee}}(r
\Lambda_0)$, then there exists a unique one-dimensional submodule in
$W$.
\end{lem}
\proof
The proof is by case by case consideration.
\begin{itemize}
\item For type ${\tt A}_n$ we have $\theta^{\vee} = \om_1^{\vee} +
\om_n^{\vee}$, so by Theorem 2
$$
\overline{V}_{-\theta^{\vee}}(r \Lam_0) \simeq V( r \om_1^{*})
\otimes V( r \om_n^{*})
$$
contains an unique one-dimensional submodule.
\item For type ${\tt B}_n$ and ${\tt D}_n$, $\theta^{\vee} =
\om_2^{\vee}$. By Theorem 2  $\overline{V}_{-w_2^{\vee}}(r \Lam_0)$
contains a unique one-dimensional submodule.
\item For type ${\tt C}_n$, $\theta^{\vee} = \om_1^{\vee}$ and
$\theta = 2 \om_1$, so again by Theorem 2
$\overline{V}_{-\theta^{\vee}}(r \Lam_0)$
contains a unique one-dimensional submodule.
\item If $\Lg$ is of type ${\tt E_6}$, ${\tt E}_7$, ${\tt E}_8$, ${\tt F}_4$,
${\tt G}_2$, then $\theta^{\vee} = \om_2^{\vee},
\om_1^{\vee},\om_8^{\vee},\om_1^{\vee},\om_2^{\vee}$ respectively and
the claim follows again by Theorem 2.
\endpf
\end{itemize}
We come to the proof of the theorem:\\
\proof Let W be the $\Lg$-module $\overline{V}_{s_{\theta}s_0}(r
\Lam_0)$. Consider the following sequence of Weyl group elements:
$$
w_0 < s_{\theta}s_0 w_0 < (s_{\theta}s_0)^{2} w_0 <
(s_{\theta}s_0)^{3}w_0 < \ldots \, .
$$
Note that the length is additive (recall $t_{-\theta} = s_{\theta}s_0$ and Lemma
\ref{cogewicht+w0}), and in a reduced decomposition of $s_{\theta}$ every
simple reflection $s_i, i=1,\ldots,n$, has to occur. So given an arbitrary element $\kappa \in \waff$,
there exists an $N\in \bn$ such that $w\le (s_{\theta}s_0)^{N}w_0$. Hence:
$$
V( \Lambda) = \lim_{ N \to \infty}
V_{(s_{\theta}s_0)^{N}w_0}(\Lambda)
$$
Write $ \Lambda = r\Lambda_0 + \lambda$, then we obtain (using the
Demazure operator)
$$
\begin{array}{rcl}
D_{(s_{\theta}s_0)^{N}w_0}( e^{(r\Lam_0 + \lambda)}) &=&
D_{(s_{\theta}s_0)^{N}}
D_{w_0}( e^{(r\Lam_0 + \lambda)}) \\
& =& D_{(s_{\theta}s_0)^{N}} ( e^{r \Lam_0} \mbox{ Char }V(\lambda) )\\
& =& e^{r \Lambda_0} ( \mbox{ Char }W)^{N} \mbox{ Char }V(\lambda)
\end{array}
$$
This shows that in the sequence of inclusions
$$
V(\lambda) \hookrightarrow W \otimes V(\lambda) \hookrightarrow W
\otimes W \otimes V(\lambda) \hookrightarrow \ldots
$$
the submodules $W^{\otimes N} \otimes V(\lambda)$ are, as
$\Lg$--modules, isomorphic to $\overline{V}_{(s_{\theta}s_0)^{N}
w_0}(\Lambda)$. Now the same arguments as in \cite{Magyar},
chapter 3, prove the theorem.
\endpf

\section{The twisted case}
In this section we would like to extend the results to twisted
affine Kac-Moody algebras and by the way to so called special
vertices. Let $\tt{X}_{n}^{(r)}$ be Dynkin diagram of affine type,
$r$ the order of the automorphism, in this section we consider
only $ r >1$. A vertex $k$ of the Dynkin diagram is called special
if $\delta - a_k \alpha_k$ is a positive root, here $\delta, a_k,
\alpha_k$ and so on are defined in the same way as in chapter 2.
For example, $0$ is always special vertex, one has $a_0 = 2$ for
$\tt {A}_{2l}^{(2)}$ and $a_0 = 1$ for the other case.

Suppose $k$ is a special vertex. Set $\theta_k = \delta - a_k
\alpha_k$, we have the finite Weyl group $W_k = \langle s_i \, | \,
i \neq k \rangle$ and let $M_k$ be the $\mathbb{Z}$-lattice spanned
by $\nu(W_k (\theta_k^{\vee}))$ (see \cite{Kac} for more details).
One knows (\cite{Kac}) that the affine Weyl group of
$\tt{X}_{n}^{(r)}$ is isomorphic to $W_k \ltimes t_{M_k}$, the
semi-direct product of $W_k$ with the translations (modulo $\delta$)
by $M_k$. The following Lemma holds:
\begin{lem}\label{specialvertex}
Let $k$ be a special vertex, then $s_k s_{\theta_k} =
t_{\nu(\beta^{\vee})}$ modulo $\delta$. For $\lam$ with $\langle
\lam, K\rangle = 0$ it follows:
$$
s_{k} s_{\theta_k}(\lam) = \lam
$$
\end{lem}
\proof
$$
\begin{array}{rcl}
s_{k}s_{\theta_k} (\lam) & = &s_{k} ( \lam - \lam(\theta_k^{\vee})\theta_k)\\
& = & s_{k}( \lam - \lam(\theta_k^{\vee})(\delta - a_k \alpha_k))\\
& = & \lam - (\lam(\alpha_k^{\vee}) + a_k \lam(\theta_k^{\vee})) \alpha_k - \lam(\theta_k^{\vee})\delta\\
\end{array}
$$
So the lemma follows, because $\lam(\alpha_k^{\vee} +
a_k\theta_k^{\vee}) = 0$, since $\lam(K) =0$.
\endpf
In section 2 we have defined the Demazure operator $D_{\beta}$ for
every real root $\beta$, with Lemma \ref{demazurezerleg} and
Lemma \ref{specialvertex} it follows:
\begin{lem}
Let $\chi \in \mathbb{Z}[\ph\cap\Lh^*]$. If $s_{\theta_k}(\chi)=
\chi$, then $D_{\al_k}(\chi)=\chi$.
\end{lem}
If one deletes in the Dynkin diagram of $\tt{X}_{n}^{(r)}$ the zero
node,  then one gets the diagram (let us call it $\tt{Y}_{n}$) of a
simple Lie Algebra. The following list shows which diagram one gets
after removing the zero node, and further, it shows that the
positive root $\delta - a_0 \alpha_0$ is a root of $\tt{Y}_{n}$.
\begin{itemize}
\item for $\tt{A}_{2}^{(2)}$: $\tt{A}_{1}$ and $\delta - a_0
\alpha_0 = \alpha_1 = \theta$, the highest root of $\tt{A}_{1}$
\item for $\tt{A}_{2l}^{(2)}$: $\tt{C}_{l}$ and $\delta - a_0
\alpha_0 =  \theta^l$, the highest long root of $\tt{C}_{l}$%
\item for $\tt{A}_{2l-1}^{(2)}$: $\tt{C}_{l}$ and $\delta - a_0
\alpha_0 =  \theta^s$, the highest short root of
$\tt{C}_{l}$%
\item for $\tt{D}_{l+1}^{(2)}$: $\tt{B}_{l}$ and $\delta - a_0
\alpha_0 = \theta^s$, the highest short root of $\tt{B}_{l}$%
\item for $\tt{E}_{6}^{(2)}$: $\tt{F}_{4}$ and $\delta - a_0
\alpha_0 =  \theta^s$, the highest short root of $\tt{F}_{4}$%
\item for $\tt{D}_{4}^{(3)}$: $\tt{G}_{2}$ and $\delta - a_0
\alpha_0 =  \theta^s$, the highest short root of $\tt{G}_{2}$%
\end{itemize}
More generally,  a vertex $k$ is special if and only if there exists
an automorphism $\sigma$ of the Dynkin diagram, such that $\sigma(k)
= 0$. In the untwisted case special is the same as minuscule. In the
twisted case, there are only for $\tt{A}_{2l-1}^{(2)}$ and
$\tt{D}_{l+1}^{(2)}$ nontrivial automorphisms. We make a new list
now for the twisted case, we delete a special vertex $k\neq 0$.
\begin{itemize}
\item for $\tt{A}_{2l-1}^{(2)}$ deleting 1: $\tt{C}_{l}$ and
$\delta - a_1 \alpha_1 =  \theta^s_1$, the highest short root of
$\tt{C}_{l}$%
\item for $\tt{D}_{l+1}^{(2)}$ deleting l: $\tt{B}_{l}$ and
$\delta - a_l \alpha_l = \theta^s_l$, the highest short root of $\tt{B}_{l}$%
\end{itemize}
We get an analog of Lemma \ref{sechs}. Let $\Lhg$
be the affine Kac-Moody algebra associated to $\tt{X}_{n}^{(r)}$,
let $\La$ be the simple Lie algebra associated to $\tt{Y}_{n}$ and
denote $P$ the weight lattice of $\La$.
\begin{lem}\label{twisted_idea}
Let $V$ be a finite dimensional $\La$ module such that $\charc V \in
\mathbb{Z}[P]$, then
\begin{equation}
D_i(\charc V) = \charc V \; \forall \, i = 0 ,\ldots ,n
\end{equation}
\end{lem}
\proof  $\charc V$ is stable under $D_i$ , $i \geq 1$. In fact,
$\charc V$ is stable under $D_\beta$ for all roots of the Lie
algebra $\La$. So only the case $i = 0$ has to be considered. Now
all weights in $V$ are of level 0, so $D_0=D_{-a_0\theta_ 0}$
$\theta_ 0 = \delta - a_0 \alpha_0$, on these weights, which
finishes the proof. This suffices to prove this, because if $\chi$
is stable under $D_{\beta}$, then it is stable under $D_{n\beta}$,
even if it is not a root.
\endpf
Recall $P$ is the $\mathbb{Z}$-lattice spanned by the fundamental
weights of $\La$. One can now formulate a statement analogous to
Theorem 1. Let $\lambda^{\vee}$ be a dominant element of $M_k
\subset P_k$, where $P_k$ are the integral, dominant weights of
$\La$. Let $\lambda^{\vee} = \lambda^{\vee}_1 + \lambda^{\vee}_2 +
\ldots + \lambda^{\vee}_r$ be a decomposition of $\lambda^{\vee}$
as a sum of dominant elements of $M_k$.
\begin{thm} Let $k$ be a special vertex of a twisted affine Kac--Moody algebra of type $\tt{X}_{n}^{(r)}$,
and let $\La,\ldots$ be as above. For all $m\geq 1$, we have an
isomorphism of $\La$-modules between the Demazure module
$\overline{V}_{-\lambda^{\vee}}(m\Lam_k)$ and the tensor product of
Demazure modules:
$$
\overline{V}_{{-\lambda^\vee}}(m\Lam_k)\simeq
\overline{V}_{{-\lambda_1^\vee}}(m\Lam_k)\otimes
\overline{V}_{{-\lambda_2^\vee}}(m\Lam_k) \otimes\cdots\otimes
\overline{V}_{{-\lambda_r^\vee}}(m\Lam_k).
$$
\end{thm}
 With the Lemma above, the proof is the same as in
the untwisted case.\\
As in the untwisted case we can now look in more detail at the
smallest Demazure modules $V_{-\om_i}(l\Lam_0)$, where $\om_i$ is
a fundamental weight for $\La$. The decompositions listed below
have been partially calculated (or conjectured) in \cite{HKOTT},
the remaining cases (and the proofs of the conjectured
decompositions) have been calculated by Naito and Sagaki
(unpublished result) as in \cite{NaitoSagaki2}
With a bar we denote again the $\La$-module,
where  $\La$ denotes the simple Lie algebra associated to diagram obtained after
removing the zero node.
Let $\epsilon = 1$ for $i$ odd and $0$ for $i$ even.
\begin{itemize}
\item $\tt{A}_{2n}^{(2)}$, $\La$ is of type $\tt{C}_{n}$
\begin{itemize}
\item[] $\overline{V_{-\om_i}(l\Lam_0)} \simeq \bigoplus\limits_{s_1 + \ldots + s_i \leq l} V(s_1 \om_1 + \ldots + s_i \om_i)$ %
\end{itemize}
\item $\tt{A}_{2n-1}^{(2)}$,$\La$ is of type $\tt{C}_{n}$
\begin{itemize}
\item[] $\overline{V_{-\om_i}(l\Lam_0)} \simeq \bigoplus\limits_{s_{p_i} + s_{p_i+2} + \ldots + s_i  = l} V(s_{p_i}\om_{p_i} + s_{p_i+2}\om_{p_i+2} +  \ldots + s_i \om_i)$
\end{itemize}
\item $\tt{D}_{n+1}^{(2)}$, $\La$ is of type $\tt{B}_{n}$
\begin{itemize}
\item[] $i= n: \overline{V_{-\om_i}(l\Lam_0)} \simeq V(l\om_n)$%
\item[] $i \neq n: \overline{V_{-\om_i}(l\Lam_0)} \simeq \bigoplus\limits_{s_1 + \ldots + s_i \leq l} V(s_1 \om_1 + \ldots + s_i \om_i)$ %
\end{itemize}
\item $\tt{E}_{6}^{(2)}$, $\La$ is of type $\tt{F}_{4}$
\begin{itemize}
\item[] $i =1: \overline{V_{-\om_i}(l\Lam_0)} \simeq \bigoplus\limits_{0 \leq s \leq l} V(s \om_1)$ %
\item[] $i =4:  \overline{V_{-\om_i}(l\Lam_0)} \simeq \bigoplus\limits_{0 \leq s_1 + s_4 \leq l} V(s_1 \om_1 + s_4\om_4)$\\%
\end{itemize}
\item $\tt{D}_{4}^{(3)}$, $\La$ is of type $\tt{G}_{2}$
\begin{itemize}
\item[] $i=1: \overline{V_{-\om_i}(l\Lam_0)} \simeq \bigoplus\limits_{0 \leq s \leq l} V(s \om_1)$ %
\end{itemize}
\end{itemize}
For the other special vertices the decompositions can be computed by taking automorphisms.\\
Theorem 3 holds in the same way, for the basic module $W$ of the
direct limit one choose $\overline{V}_{- \theta_k^{\vee}}(r
\Lam_k)$. Then the direct sum decomposition of $W$ contains
obviously an one dimensional module, namely the one who corresponds
in the Demazure module $V_{-\theta_k^{\vee}}(r \Lam_k)$ to the
weight $r \Lam_k$. Again let $V^{\infty}_{\lam,r}$ be the direct
limit constructed above. Then it follows
\begin{thm}
For any integral dominant weight $\Lam$ of $\Lhg$, $\Lam = r \Lam_k
+ \lam$, the $\La$-modules $V^{\infty}_{\lam,r}$ and $\overline{V}(r
\Lam_k)$ are isomorphic.
\end{thm}

\end{document}